\documentclass[12pt, reqno]{amsart}
\usepackage[arrow,matrix,curve]{xy}

\newtheorem{theorem}{Theorem}[section]
\newtheorem{proposition}[theorem]{Proposition}
\newtheorem{lemma}[theorem]{Lemma}
\newtheorem{corollary}[theorem]{Corollary}

\theoremstyle{definition} 
\newtheorem{definition}[theorem]{Definition}

\newtheorem{remark}[theorem]{Remark}

 \newcommand{\N}{\mathbb N}
\newcommand{\R}{\mathbb R} \newcommand{\T}{\mathbb T}
\newcommand{\Z}{\mathbb Z}

\DeclareMathOperator{\XX}{{\overline{\it X}}}
\DeclareMathOperator{\fx}{{[\overline{\it F}_{\it \! X}]}}

\DeclareMathOperator{\stsys}{{\rm stsys}}
\DeclareMathOperator{\confsys}{{\rm confsys}}
\DeclareMathOperator{\pisys}{{\rm sys}\pi}

\DeclareMathOperator{\vol}{{\rm vol}}

\DeclareMathOperator{\AJ}{{\mathcal A}} 
\DeclareMathOperator{\baraj}{{\overline{\mathcal A}}}

\DeclareMathOperator{\jac}{{\rm Jac}}

\DeclareMathOperator{\trace}{{\rm trace}}

\DeclareMathOperator{\cuplength}{{\rm cuplength}}

\def\st{{\rm st}}
\def\dR{{\rm dR}}

\def\conf{{\rm conf}}

\def\stp{{p}}

\def\dimtarget{{b}}

\def\ie {{\it i.e.\ }} 
\def\eg {{\it e.g.\ }} 
\def\cf {\hbox{\it cf.\ }}

 \def\Pr {{\rm Pr}}

\numberwithin{equation}{section} 
\numberwithin{figure}{section}

\begin{document}

\author[V.~Bangert]{Victor Bangert$^!$} \address{ Mathematisches
Institut, Universit\"at Freiburg, Eckerstr.~1, 79104 Freiburg,
Germany} \email{bangert@mathematik.uni-freiburg.de}
\thanks{$^!$Partially Supported by DFG-Forschergruppe `Nonlinear
Partial Differential Equations: Theoretical and Numerical Analysis'}

\author[C.~Croke]{Christopher Croke$^+$} \address{ Department of
Mathematics, University of Pennsylvania, Philadelphia, PA 19104-6395
USA} \email{ccroke@math.upenn.edu} \thanks{$^+$Supported by NSF grant
DMS 02-02536 and the Max-Planck-Institut f\"ur Mathematik Bonn}

\author[S.~Ivanov]{Sergei V. Ivanov$^{\dagger}$} \address{Steklov
Math. Institute, Fontanka 27, RU-191011 St. Petersburg, Russia}
\email{svivanov@pdmi.ras.ru} \thanks{$^\dagger$Supported by grants
CRDF RM1-2381-ST-02, RFBR 02-01-00090, and NS-1914.2003.1}

\author[M.~Katz]{Mikhail G. Katz$^{*}$} \address{Department of
Mathematics, Bar Ilan University, Ramat Gan 52900 Israel}
\email{katzmik@math.biu.ac.il} \thanks{$^{*}$Supported by the Israel
Science Foundation (grants no.\ 620/00-10.0 and 84/03)}

\title[Case of equality in Loewner-type inequalities] {Boundary case
of equality in optimal Loewner-type inequalities$^{1}$
}



\begin{abstract}
We prove certain optimal systolic inequalities for a closed Riemannian
manifold $(X,g)$, depending on a pair of parameters, $n$ and $b$.
Here $n$ is the dimension of $X$, while $b$ is its first Betti number.
The proof of the inequalities involves constructing Abel-Jacobi maps
from $X$ to its Jacobi torus $\T^b$, which are area-decreasing (on
$b$-dimensional areas), with respect to suitable norms.  These norms
are the stable norm of $g$, the conformally invariant norm, as well as
other $L^p$-norms.  Here we exploit $L^p$-minimizing differential
1-forms in cohomology classes.  We characterize the case of equality
in our optimal inequalities, in terms of the criticality of the
lattice of deck transformations of~$\T^b$, while the Abel-Jacobi map
is a harmonic Riemannian submersion.  That the resulting inequalities
are actually nonvacuous follows from an isoperimetric inequality of
Federer and Fleming, under the assumption of the nonvanishing of the
homology class of the lift of the typical fiber of the Abel-Jacobi map
to the maximal free abelian cover.
\end{abstract}

\footnotetext[1]{\large {\em Transactions Amer. Math. Soc.}, to
appear}

\maketitle

\tableofcontents

\section{Introduction, conjectures, and some results}

In the study of optimal systolic inequalities on a Riemannian manifold
$(X,g)$, generalizing the Loewner inequality, two questions arise:
\begin{enumerate}
\item
Which stable systolic inequalities can be replaced by their conformal
analogues?
\item
Calculate the values of the constants arising in the inequalities.
\end{enumerate}
The definitions of the systolic invariants $\stsys_1(g),
\confsys_1(g)$, and $\pisys_1(g)$ can be found in the survey
\cite{CK}, while the Abel-Jacobi map $\AJ_X$, in \cite{Li, Gr2}, \cf
\cite[(4.3)]{BaKa03}, \cite{IK}.  The first such optimal inequality
was proved in~\cite [pp.~259-260] {Gr3} (\cf \cite[inequality
(5.14)]{CK}) based on the techniques of D. Burago and the third author
\cite{BI94, BI95}.
\begin{theorem}[M. Gromov]
\label{11a}
Let $X$ be a closed orientable manifold.  Let $n=\dim(X)$ and assume
$n=b_1(X)=\cuplength_\R(X)$.  Then every metric $g$ on $X$ satisfies
the following optimal inequality:
\begin{equation}
\label{10}
\stsys_1(g)^n \leq \deg(\AJ_X)^{-1} \left( \gamma_n
\right)^{\frac{n}{2}} \vol_n(g),
\end{equation}
where $\gamma_n$ is the Hermite constant of \eqref{34}, while $\AJ_X:
X \to \T^n$ is the classifying map of the natural epimorphism
$\pi_1(X)\to H_1(X,\Z)/{\rm torsion}$.
\end{theorem}
The following conformally invariant generalisation of Gromov's
inequality \eqref{10} is immediate from Proposition~\ref{21a}, \cf
inequality~\eqref{28}.
\begin{theorem}
\label{11c}
Let $X$ be a closed orientable manifold.  Let $n=\dim(X)$ and assume
$n=b_1(X)=\cuplength_\R(X)$.  Then every metric $g$ on $X$ satisfies
the following optimal inequality:
\begin{equation}
\label{10c}
\confsys_1(g)^n \leq \deg(\AJ_X)^{-1} \left( \gamma_n
\right)^{\frac{n}{2}} .
\end{equation}
\end{theorem}

In connection with question 1 above, we aim for the following results.
We show in section \ref{vb} that if the $L^p$ norm is equal to the
comass norm, then the minimizing measurable 1-form is in fact
continuous and has constant pointwise norm.  Then we can apply a
Voronoi-type argument to prove that in an extremal situation, in fact
all cohomology classes are represented by 1-forms of constant norm.
This allows us to prove that a manifold $X$ such that $\dim(X) =
b_1(X) = \cuplength_\R(X)$, which satisfies the boundary case of
equality of (Gromov's) inequality relating the stable 1-systole and
the volume, must be a critical flat torus (Theorem~\ref{35}).  The
proof passes via the conformally invariant strengthening of Gromov's
inequality, namely Theorem \ref{11c}.

Given a Euclidean lattice $L\subset \R^n$, denote by $\lambda_1(L)>0$
the least length of a nonzero vector in $L$.  A lattice $L$ is called
{\em critical\/} if it realizes the supremum
\begin{equation}
\label{34}
\sup_{L\subset \R^n}\frac{\lambda_1(L)^n}{\vol(\R^n/L)} =
(\gamma_n)^{\frac{n}{2}},
\end{equation}
where $\gamma_n$ is the Hermite constant.  Equivalently, $L$ is
critical if the $L$-centered packing by balls is the densest one in
dimension $n$.  Recall that a critical (or, more generally, extremal)
lattice is necessarily perfect and eutactic \cite{Bar}.  The following
two theorems are proved in Section~\ref{six}.
\begin{theorem}
\label{13a}
The boundary case of equality in the conformally invariant inequality
\eqref{10c}, occurs if and only if $X$ is a torus and the metric $g$
is conformal to a flat critical metric.
\end{theorem}
\begin{theorem}
\label{35}
In the hypotheses of Theorem~\ref{11a}, equality in~\eqref{10} is
attained precisely when $X$ is a torus, while $g$ is a flat metric
whose deck transformations form a critical lattice.
\end{theorem}

In Section~\ref{smallcycles}, we discuss a notion of degree of the
Abel-Jacobi map when the dimension $n$ of the manifold $X$ exceeds its
first Betti number~$b$.  Here we prove the positivity of the degree
under a suitable topological hypothesis of homological nonvanishing.
In Section~\ref {thre}, we extend Theorem~\ref{34} and
Theorem~\ref{35} to a general pair $n\geq b$, and discuss the related
results in the literature.  In Section~\ref{three}, we define the
norms involved, and state a key proposition used in the proof of
Theorems \ref{34} and \ref{35}.  In Sections~\ref{vb} and
\ref{eleven}, we collect analytical results pertaining to closed forms
that minimize an $L^p$-norm.  This constitutes the analytical backbone
of the paper.  Sections \ref{bi} and \ref{eight} contain successive
generalisations of a construction of \cite{BI94, BI95}.  The remaining
sections are devoted to the proof of the results of
Section~\ref{thre}.

\section{Notion of degree when dimension exceeds Betti number}
\label{smallcycles}

Let $(X,g)$ be a closed Riemannian manifold, $b=b_1(X)$.  We denote by
$\AJ_X$ the Abel-Jacobi map,
$$
  \AJ_X:X\to H_1(X,\R)/H_1(X,\Z)_\R \simeq \T^b ,
$$
inducing (the natural) isomorphism in 1-dimensional cohomology.
Let~$\XX$ be the maximal free abelian cover of $X$, \cf
Definition~\ref{xx}.  Denote by~$\baraj_X$ the proper map obtained as
the lift of $\AJ_X$ to $\XX$:
$$
\baraj_X:\XX\to H_1(X,\R) \simeq \R^b .
$$
If $\XX$ is orientable we fix an orientation for $\XX$.  In the case
that $X$ is orientable, the following definition is due to
M. Gromov~\cite[p.~101]{Gr1}.
\begin{definition}
\label{deg}
Let $\fx\in H_{n-b}(\XX,R)$ denote the homology class of the regular
fibers $\baraj_X^{-1}(y)$ of $\baraj_X$ where we take $R=\Z$ if $\XX$
is orientable and $R=\Z_2$ otherwise.  Then the geometric degree
\[
\deg(\AJ_X)
\]
of $\baraj$ is the infimum of the $(n-b)$-volumes of all cycles
representing~$\fx$.
\end{definition}

\begin{remark}
If $\XX$ is orientable, then changing the orientation of $\XX$ will
not change $\deg(\AJ_X)$.  If $n=b$ and $X$ is orientable and
connected, then $\deg(\AJ_X)$ is (the absolute value of) the topological
degree of $\AJ_X$.  If $n>b$ and $\fx \not=0$, then $\deg(\AJ_X)$
depends on the metric on $X$.
\end{remark}

For our arguments the following obvious consequence of
Definition~\ref{deg} will be crucial.

\begin{lemma}
Suppose $f: X \to H_1(X,\R)/H_1(X,\Z)_\R$ is homotopic to~$\AJ_X$ and
$y$ is a regular value of $f$.  Then we have
\[
\vol_{n-b}(f^{-1}(y)) \geq \deg(\AJ_X).
\]
\end{lemma}

\begin{proof}
Let $\bar y \in H_1(X,\R)$ be a preimage of $y$.  Let $\bar f: \XX \to
H_1(X,\R)$ be a lift of $f$.  Then the covering projection $\XX \to X$
maps the fiber~$\bar f^{-1}(\bar y)$ isometrically onto the fiber
$f^{-1}(y)$, so that we have
\[
\vol_{n-b}(\bar f^{-1}(\bar y)) = \vol_{n-b}(f^{-1}(y)) .
\]
On the other hand, $\bar f$ is homotopic to $\baraj_X$ by a proper
map, and hence the fiber~$\bar f^{-1}(\bar y)$ represents the class
$\fx$.
\end{proof}

All the results discussed below are nonvacuous only if~$\deg(\AJ_X)$
is positive.  So it is important to know when this is the case.  First
note that this is a topological property.

\begin{remark}
If $\fx\not = 0$ then $\deg(\AJ_X)> 0$.
\end{remark}

This follows from the deformation theorem as stated in
\cite[4.2.9]{Fe}, complemented by \cite[$(4.2.9)^\nu$]{Fe} for the
case of $\Z_2$-coefficients, see also \cite{FF}, \cite{Whi}, and
\cite[Chapter~5]{Mo} for a helpful picture.  The deformation theorem
implies that for every compact Riemannian manifold~$X$, there exists
$\epsilon > 0$ such that for every $q\in \N$ every $q$-cycle
of~$q$-volume smaller than $\epsilon$ can be deformed into the
$(q-1)$-skeleton of a triangulation of $X$.  Now, if $z$ is a
$q$-cycle in $\XX$ of $q$-volume smaller than $\epsilon$, then the
projection of $z$ to $X$ can be deformed into the $(q-1)$-skeleton.
Lifting this deformation to $\XX$ we conclude that the homology class
of~$z$ vanishes regardless of the ring of coefficients.

If $X$ is the quotient of the 3-dimensional Heisenberg group by its
subgroup of integral elements, then the fibers of $\AJ_X$ are
homologically trivial in $X$ while $\fx\not=0$.  This explains why one
looks at the situation in $\XX$ in Definition~\ref{deg}.

In analogy with the maximal $\cuplength_\R$ condition in
Theorem~\ref{11a}, we note that $\fx\not=0$ if $X$ is orientable and
if $\alpha_1\cup\ldots \cup \alpha_b \not= 0$ for every basis
$\alpha_1,\ldots, \alpha_b$ of $H^1(X,\R)$.  Even if the real cup
product vanishes, for a 3-manifold $X$ with Betti number 2, the
nonvanishing of the self-linking number of the typical fiber
$\AJ_X^{-1}(p)$ of $\AJ_X: X \to \T^2$ is a sufficient condition for
the nonvanishing of the fiber class $\fx$, see~\cite{KL}.  Note that,
if the minimizing rectifiable current in $\fx$ happens to be
nonsingular, then the work \cite{Lcirce} produces explicit lower
bounds for~$\deg(\AJ_X)$ in terms of the injectivity radius and an
upper bound on the sectional curvature.

\section{The theorems}		
\label{thre}

The following two results were proved in \cite{IK}.
\begin{theorem}
\label{12}
Let $X$ be a closed Riemannian manifold.  Let $n=\dim(X)$ and
$b=b_1(X)$, and assume $n\ge b\geq 1$.  Then every metric $g$ on~$X$
satisfies the inequality
\begin{equation}
\label{eq12}
{\deg(\AJ_X) \stsys_1(g)^b} \leq (\gamma_b) ^{\frac{b}{2}} {\vol_n
(g)} .
\end{equation}
\end{theorem}

\begin{corollary}
\label{13}
Let $X$ be a closed orientable manifold.  Let $b=b_1(X)$.  Assume that
$\dim(X)=b+1$, and $\fx \not=0$.  Then every metric $g$ on~$X$
satisfies the following optimal inequality:
\begin{equation}
\label{11}
\stsys_1(g)^b \pisys_1(g) \leq (\gamma_b)^{\frac{b}{2}} \vol_{b+1}(g).
\end{equation}
\end{corollary}

The main tool in studying the case of equality of the above
inequalities is Proposition~\ref{23b} which concerns a kind of an
$L^p$-systole, \cf Definition~\ref{53}.  Let $p\geq 1$.

\begin{definition}
\label{23z}
The norm $\|\;\|_p$ on $H_1(X;\R)$ is the dual norm to the~$L^p$-norm
$\|\;\|^*_p$ on $H^1(X;\R)$, \cf \cite{BK}.
\end{definition}

\begin{remark}
\label{23d}
In general, if $h\in H_k(X,\R)$, then the quantity $\| h
\|_{\frac{n}{k}}$ is the supremum of the stable norms of $h$ with
respect to all metrics that are conformal to $g$ and have unit volume,
\cf \cite [7.4.A] {Gr1}.  In particular, one has
\begin{equation}
\label{28}
\stsys_k(g) \le \confsys_k(g) \vol_n(g)^{\frac{k}{n}}.
\end{equation}
\end{remark}

Given a lattice $L$ equipped with a norm $\|\;\|$, denote by
$\lambda_1(L,\|\;\|)$ the least norm of a nonzero vector in $L$.  The
following proposition generalizes Theorem~\ref{12}.
\begin{proposition}
\label{23b}
Let $X$ be a closed orientable Riemannian manifold of unit volume.
Let $b=b_1(X)$.  Let $p \geq \max\{b, 2 \}$.  Then the following
inequality is satisfied:
\begin{equation}
\label{23}
\deg(\AJ_X) \; \lambda_1\! \left( H_1(X,\Z)_\R ^{\phantom{I}},
\|\;\|_p \right) ^b \leq (\gamma_b) ^{\frac{b}{2}}.
\end{equation}
\end{proposition}
The proof of Proposition~\ref{23b} appears in Section~\ref{thirteen}.
By choosing~$p=n$ in Proposition~\ref{23b}, we obtain the following
generalisation of Theorem~\ref{11c} and Theorem~\ref{12}.

\begin{theorem}
\label{37}
Let $X$ be a closed orientable manifold.  Let $n=\dim(X)$ and
$b=b_1(X)$.  Then every metric $g$ on $X$ satisfies the inequality
\begin{equation}
\label{23c}
\deg(\AJ_X) \confsys_1(g) ^b \leq (\gamma_b) ^{\frac{b}{2}}
\vol_n(g)^{n-b}.
\end{equation}
\end{theorem}

\begin{remark}
In the case $b=2$, Theorem~\ref{37} follows by setting $p=2$, and
applying the coarea formula to Lichnerowicz's harmonic map~\cite{Li},
see Remark~\ref{121}.  This proves inequality~\eqref{23} for the
``$L^2$-systole'', and the monotonicity in $p$ then proves
\eqref{23c}.  One can thus avoid in this case the extensions of the BI
construction, described in Section~\ref{bi} and Section~\ref{eight}.
Thus our constructions can be viewed as generalisations, to~$b\geq 3$,
of the the properties of being area-decreasing (on the average) of
Lichnerowicz's map in the case $b=2$.
\end{remark}
In analogy with our Theorem~\ref{35}, we prove the following.
\begin{theorem}
\label{17}
Let $X$ be a closed orientable manifold.  Assume that
$\dim(X)>b_1(X)$.  Suppose equality is attained in~\eqref{23c}.  Then
\begin{enumerate}
\item
all harmonic one-forms on $X$ have pointwise constant norm;
\item
$X$ admits a harmonic Riemannian submersion onto a flat $b$-torus
whose deck group is a critical lattice in $\R^b$;
\item
the submersion is given by the Abel-Jacobi map defined via harmonic
one-forms;
\item
the fibers of the submersion are minimal surfaces of constant
$(n-b)$-dimensional area.
\end{enumerate}
\end{theorem}

Note that Riemannian manifolds for which the Abel-Jacobi map is a
harmonic Riemannian submersion onto a flat torus are studied
in~\cite{BaKa03}.  Harmonic Riemannian submersions to flat tori are
examples of harmonic morphisms, see \eg the monograph \cite{BW}.

As a consequence of Theorem~\ref{17}, we obtain the following, see
some related results in~\cite{NV, Na}.
\begin{corollary}
In the hypotheses of Corollary~\ref{13}, equality in \eqref{11} is
attained precisely when $X$ is a two-step nilpotent manifold with
1-dimensional center, while $g$ is the metric of a Riemannian
submersion with geodesic fibers of constant length.
\end{corollary}

The case $n=b+2$ is studied in \cite{IK} and the companion
paper~\cite{BCIK1}, which explores the closely related filling area
conjecture.  A different higher-dimensional generalisation of the
Loewner inequality is studied in \cite{BK, BaKa03}.  A proof of the
Loewner inequality in genus~2 appears in \cite{KS1}, while the
asymptotics of the optimal systolic ratio for large genus are studied
in \cite{KS2}.  Similar asymptotics for the conformal 2-systole are
studied in \cite{[K2]}.  The work \cite{KR} provides a general
framework for systolic geometry, in terms of a notion of {\em systolic
category}, related to Lusternik-Schnirelmann category.

\section{Stable norms and conformal norms}
\label{three}

Let $X$ and $Y$ be compact smooth manifolds.  Let $\varphi:X\to Y$ be
a continuous map inducing an epimorphism in one-dimensional real
homology.  Given a Riemannian metric on $X$ (more generally, the
structure of a locally simply connected length space), one defines the
{\it relative stable norm} $\|\;\|_{\st/\varphi}$ on $H_1(Y;\R)$ by
setting
\begin{equation}
\|\alpha\|_{\st/\varphi} = \inf \left\{ \left. \|\beta\|_{\st}
^{\phantom{I}} \; \right| \beta\in H_1(X;\R),\ \varphi_*(\beta)=\alpha
\right\} ,
\end{equation}
where $\|\;\|_{\st}$ is the ordinary (``absolute'') stable norm of
$X$.  The stable norm itself may be thought of as the relative stable
norm defined by the Abel-Jacobi map to the torus $H_1(X,\R)/ H_1 (X,
\Z) _{\R}$.

Let $X^n$ be a compact Riemannian manifold, $V^\dimtarget$ a vector
space, and~$\Gamma$ a lattice in $V$.  We will identify $V$ and
$H_1(V/\Gamma;\R)$.  Let $\varphi:X\to V/\Gamma$ be a continuous map
inducing an epimorphism of the fundamental groups.  Note that there is
a natural isomorphism 
\[
\Gamma\simeq\pi_1(V/ \Gamma) \simeq \pi_1(X)/ker(\varphi_*).
\]
The following proposition was proved in \cite{IK}.  Denote by
$\|\;\|_E$ the Euclidean norm on $V$, which is defined by the John
ellipsoid, \cf \cite{Jo}, of the relative stable norm~$\| \; \| _{
\st/ \varphi}$.

\begin{proposition} \label{area-nonexp}
There exists a Lipschitz map $X\to(V/\Gamma, \|\; \| _E)$ which is
homotopic to~$\varphi$ and non-expanding on all
$\dimtarget$-dimensional areas, where $\dimtarget =\dim V$.
\end{proposition}

The following statement may be known by convex set theorists.  Anyway,
there is a proof in \cite{BI94}.

\begin{lemma} 
\label{ellipsoid-decomposition}
Let $(V^\dimtarget , \|\;\|)$ be a Banach space.  Let $\|\;\|_E$ be
the Euclidean norm determined by the John ellipsoid of the unit ball
of the norm~$\|\;\|$. Then there exists a decomposition of
$\|\;\|_E^2$ into rank-1 quadratic forms:
$$
  \|\;\|_E^2 = \sum_{i=1}^N \lambda_i L_i^2
$$
such that $N\le \frac{\dimtarget (\dimtarget +1)}{2}$, $\lambda_i>0$
for all $i$, $\sum\lambda_i= \dimtarget$, and $L_i:V\to\R$ are linear
functions with $\|L_i\|^*=1$ where $\|\;\|^*$ is the dual norm
to~$\|\;\|$.
\end{lemma}

Note that the bound $N\le \frac{\dimtarget (\dimtarget +1)}{2} +1$
appears in \cite{BI94}, but one easily sees that it can be reduced by
one.

Given a Riemannian manifold $X^n$, one defines a conformally invariant
norm $\|\;\|_{\conf}$ on $H_1(X;\R)$ as the dual norm to the
$L^n$-norm~$\|\;\|^*_n$ on~$H^1(X;\R)$, \cf Remark \ref{23d}.  If
$\vol(X)=1$, then we have an inequality $\|\;\|^*_n \le \|\;\|^*$,
where $\|\;\|^*$ is the comass norm in cohomology, dual to the stable
norm in homology.  Therefore dually, we have the
inequality~$\|\;\|_{\conf}\ge \|\;\|_{\st}$, \cf inequality
\eqref{28}.  In particular, the volume form defined by the John
ellipsoid of the conformal norm is greater than or equal to that of
the stable norm.

For a continuous map $\varphi:X\to Y$ inducing an epimorphism in
one-dimensional homology, one defines a relative conformal norm
$\|\;\|_{\conf/\varphi}$ on $H_1(Y;\R)$ by setting
\begin{equation}
\|\alpha\|_{\conf/\varphi} = \inf \left\{
\left. \|\beta\|_{\conf}^{\phantom{I}} \; \right| \beta\in H_1(X;\R),\
\varphi_*(\beta)=\alpha \right\} .
\end{equation}
We present an ``integral'' version of Proposition \ref{area-nonexp}
for conformal norms.

\begin{proposition}
\label{21a}
Let $X^n$ be a compact Riemannian manifold, $V^n$ a vector space, and
let $\Gamma$ be a lattice in $V$.  Let $\varphi:X\to V/\Gamma$ be a
continuous map inducing an epimorphism of the fundamental groups.  Let
$\|\;\|_E$ denote the Euclidean norm on $V$ defined by the John
ellipsoid of the relative conformal norm $\|\;\|_{\conf/\varphi}$.
Then there exists a $C^1$ map
\[
f:X\to(V/\Gamma,\|\;\|_E)
\]
which is homotopic to $\varphi$ and satisfies
\[
\int_X \jac(f) \le 1.
\]
In particular, if $\varphi$ has nonzero degree (or nonzero absolute
degree), then we have the following upper bound for the volume of the
torus $V/\Gamma$:
\[
\vol(V/\Gamma,\|\;\|_E)\le \frac{1}{\deg(\varphi)}.
\]
\end{proposition}

An adaptation of the main construction \cite{BI94,BI95} involved in
the proof will be described in Section~\ref{bi}.  We will then use it
to prove Proposition~\ref{21a} in Section~\ref{five}.

\section{Existence of $L^p$-minimizers in cohomology classes}
\label{vb}

In this section, we discuss the existence of a (weakly closed)
$L^p$-form that minimizes the $L^p$-norm, $1<p<\infty$, in a given
cohomology class, used in the proof of Proposition~\ref{21a} in
Section~\ref{five}.

Although this is known, we could not find a good reference.  The proof
consists in a straightforward application of the direct method of the
calculus of variations.  Actually C.~Hamburger \cite{Ha92} even proved
H\"older continuity of the minimizers.

Note that the existence of a continuous 1-form minimizing the comass,
and defining a geodesic lamination, is proved in \cite [Theorem~1.7]
{FS}.

If $E\rightarrow X$ is a Euclidean vector bundle over a compact,
oriented Riemannian manifold $X$, and if $1\le p<\infty$, let $L^p(E)$
denote the vector space of $L^p$-sections $s$ of $E$ endowed with the
Banach norm
\begin{equation}
|s|_p = \left( \int_X |s(x)|^p d\vol (x) \right)^{\frac{1}{p}}.
\end{equation}
If $p=\infty$, let $L^{\infty}(E)$ denote the essentially bounded
sections of $E$ and let $|s|_{\infty}$ denote the essential supremum
of the function $x\in X\rightarrow|s(x)|$.

In analogy with the case of functions $X\rightarrow\R$, we have: 

\begin{lemma}
\label{Lemma1.1.} 
{\it If $1\le p<\infty$ and $\frac{1}{p} +\frac{1}{q}=1$, then the
natural map
$$
J_p : L^p(E)\rightarrow L^q(E)^*, (J_p(s))(s'):=\int_X\langle
s,s'\rangle d\vol
$$
for $s\in L^p(E)$, $s'\in L^q(E)$ is an isometric isomorphism between
$L^p(E)$ and the dual space $L^q(E)^*$ of the Banach space $L^q(E)$.
For $1<p<\infty$ the Banach space $L^p(E)$ is reflexive. }
\end{lemma}

\begin{proof}
If we add a vector bundle $E'$ to $E$ such that $E\oplus E'\rightarrow
X$ is trivial, we can easily reduce the claim to the function
case. For this case see \eg \cite[4.12 and 6.10]{Al}.
\end{proof}

In the sequel, the vector bundle $E$ will be the bundle $\Lambda^k X$
of alternating $k$-forms, $0\le k\le \dim X =:n$.  The Euclidean
structure on~$\Lambda^k X$ will be the one induced by the Riemannian
metric on $X$.  For~$k\in\{0,1,n-1,n\}$ the corresponding Euclidean
norm on the fibers of $\Lambda^k X$ coincides with the comass norm,
\cf \cite[1.8.1]{Fe}.  We assume that the total volume of $X$ is
normalized to one.

Denote by $\Omega^k X$ the space of smooth $k$-forms on $X$.  Given a
cohomology class $\alpha\in H^k(X,\R)$, denote by $\alpha_{\dR}$ the
set of smooth closed $k$-forms $w\in\Omega^k X$ that represent
$\alpha$.  Then
$$
|\alpha|^*_p := \inf \left\{ \left. |w|_p^{\phantom{I}} \right|
 w\in\alpha_{\dR} \right\}
$$
defines a norm on the finite-dimensional vector space $H^k(X,\R)$.  To
see that $|\alpha|^*_p > 0$ if $\alpha\not= 0$, note that there exists
$\beta\in H^{n-k}(X,\R)$ such that the cup product $\alpha\cup\beta\in
H^n(X,\R)$ is non-zero.  This implies that
$$
\left|
\int_X w \wedge\pi
\right|
= 
\left|(\alpha\cup\beta) ([X])
\right|>0
$$
whenever $w\in\alpha_{\dR}$ and $\pi\in\beta_{\dR}$.  Because of the
pointwise inequality~$|w\wedge\pi|\le{n\choose k}^{\frac{1}{2}} |w|
|\pi|$, \cf \cite[1.7.5]{Fe} we obtain
$$
\int_X |w| \; |\pi| d \mbox{vol} \ge {n\choose k}^{-\frac{1}{2}}
|(\alpha\cup\beta) ([X])|.
$$
Hence the H\"older inequality implies
$$
|\alpha|^*_p \; |\beta|^*_q \ge {n\choose k}^{-\frac{1}{2}}
|(\alpha\cup\beta) ([X])| > 0.
$$
We let $\overline{\alpha}^p_{\dR} \subseteq L^p(\Lambda^k X)$ denote
the closure of $\alpha_{\dR}$ with respect to the~$L^p$-norm.  Note
that
\begin{equation}
\label{(1.1)}
|\alpha|^*_p = \inf \left\{ |w|_p \left| w\in {\overline {\alpha}}
^{p\phantom{I}} _{\dR} \right. \right\}
\end{equation}
Moreover, if $\beta\in H^{n-k}(X,\R), \frac{1}{p} + \frac{1}{q} = 1,
w\in \overline{\alpha}^p_{\dR}$ and $\pi\in\overline{\beta}^q_{\dR}$,
then
\begin{equation}
\label{(1.2)}
{\int_X w\wedge\pi = [\alpha\cup\beta](X).}
\end{equation}

Since the cup product $\cup:H^k(X,\R)\times H^{n-k} (X,\R) \rightarrow
\R$ is non-degenerate this implies:

\begin{lemma}
\label{Lemma1.2.} 
{\it Every smooth closed form in $\overline{\alpha} ^p_{\dR}$
represents $\alpha$.}
\end{lemma}

Next we prove that $\overline{\alpha}^p_{\dR}$ contains a unique
$L^p$-form with minimal~$p$-norm if $1<p<\infty$.

\begin{proposition}
\label{Proposition1.3.} 
{\it Let $1<p<\infty$ and $k\in\{1,\ldots,\dim X-1\}$.  For every
$\alpha\in H^k(X,\R)$ there exists a unique
$w\in\overline{\alpha}^p_{\dR}$ such that }
$$
|w|_p = |\alpha|^*_p.
$$
\end{proposition}

\begin{proof}
Let $(w_i)_{i\in\N}$ be a sequence in $\alpha_{\dR}$ such that
$\lim\limits_{i\rightarrow\infty} |w_i|_p = |\alpha|^*_p$.  Since
$L^p(\Lambda^k X)$ is a reflexive Banach space, \cf
Lemma~\ref{Lemma1.1.}, we can assume that the sequence
$(w_i)_{i\in\N}$ converges weakly to some $w\in L^p(\Lambda^k X)$, \cf
\cite[6.9]{Al}.  Since $\overline{\alpha}^p_{\dR}$ is a closed affine
subspace of $L^p(\Lambda^k X)$, the Hahn-Banach theorem implies that
$\overline{\alpha}^p_{\dR}$ is also weakly sequentially closed, \cf
\cite[6.12]{Al}.  This implies that $w\in\overline{\alpha}^p_{\dR}$.
Moreover, the weak convergence of $(w_i)_{i\in\N}$ to $w$ implies that
$$
|w|_p \le \lim\limits_{i\rightarrow\infty} |w_i|_p = |\alpha|^*_p.
$$
Using equality \eqref{(1.1)} and $w\in\overline{\alpha}^p_{\dR}$, we
conclude that $|w|_p = |\alpha|^*_p$.  The uniqueness of the minimizer
$w$ follows from the fact that the equality~$|w+\tilde{w}|_p = |w|_p +
|\tilde{w}|_p$ can hold only if $w$ and $\tilde{w}$ are linearly
dependent in $L^p(\Lambda^k X)$.~\end{proof}

\section{Existence of harmonic forms with constant norm}
\label{eleven}

We continue to assume that the volume of $X$ is normalized to one.
Then, as a direct consequence of the H\"older inequality, the
function
\[
p\in[1,\infty] \rightarrow |\alpha|^*_p
\]
is (weakly) monotonically increasing for every class $\alpha\in
H^k(X,\R)$, where $0<k< n = \dim X$.

\begin{proposition}
\label{125} 
Let $\alpha\in H^k(X,\R)$, $0<k< n$.  Assume there exist reals $p<p'$
in $[1,\infty]$ such that $|\alpha|^*_p = |\alpha|^*_{p'}$.  Then the
following holds:

\begin{enumerate}
\item
the function $q\in[1,\infty]\rightarrow |\alpha|^*_q$ is constant;
\item
the harmonic representative $\omega$ of $\alpha$ has constant
(pointwise) norm;
\item
we have $|\omega | _q = |\alpha|_q^*$ for all $q\in [1,\infty ]$.
\end{enumerate}
Conversely, if $p\in [1,\infty)$ and if there exists a form $\omega
\in \bar \alpha_{\dR}^p$ such that~$| \omega |_p = | \alpha |_p^*$ and
$\omega$ has constant norm almost everywhere, then
\[
|\omega |_p = | \omega |_q = |\alpha|_q^*
\]
for all $q\in [1,\infty]$.  In particular, $\omega$ is harmonic.
\end{proposition}

\begin{proof}
Let $p < p'$ in $[1,\infty]$ satisfy $| \alpha |_p^* = | \alpha
|_{p'}^*$.  Since $q \mapsto | \alpha |_q^*$ is monotonic, we may
assume that $1 < p < p' < \infty$.  According to Proposition~\ref
{Proposition1.3.}, there exists $\omega \in\overline{\alpha}
^{p'}_{\dR}$ such that $|\omega |_{p'}=|\alpha|^*_{p'}$.  First we show
that $\omega$ has constant norm.  Since $\overline{\alpha}^{p'}_{\dR}
\subseteq\overline {\alpha} ^{p}_{\dR}$, we see that
$$
|\alpha|^*_p \le |\omega |_p \le |\omega |_{p'} = |\alpha|^*_{p'}.
$$
Hence our assumption $| \alpha |_p^* = | \alpha |_{p'}^*$ implies
$|\omega |_p = |\omega |_{p'}$.  Therefore the function $x\in
X\rightarrow|\omega _x|$ is almost everywhere constant.

Next, we prove that $\omega$ is indeed harmonic.
If $p\leq 2$, this follows directly from 
\[
|\alpha|_p^* = | \omega |_p = |\omega |_2 \geq |\alpha |_2^* \geq
|\alpha |_p^* .
\]
For general $p$, we will now show that $\omega$ is weakly closed and
coclosed, i.e.~that
\begin{equation}
\label{(2.1)} 
{\int_X \langle \omega, d^* \pi\rangle d \,\mbox{vol} = 0 \hbox{\ for\
all\ } \pi\in\Omega^{k+1}X,}
\end{equation}
and
\begin{equation}
\label{(2.2)}
{\int_X \langle \omega ,d\sigma \rangle d \,\mbox{vol} = 0 \hbox{\
for\ all \ } \sigma\in\Omega^{k-1}X.}
\end{equation}
To prove \eqref{(2.1)}, consider a fixed $\pi\in\Omega^{k+1} X$ and
the linear functional
$$
\tilde{\omega }\in L^p(\Lambda^k X)\rightarrow\int_X\langle\tilde
{\omega}, d^* \pi \rangle d \,\mbox{vol}
$$
on $L^p(\Lambda^k X)$. This functional vanishes on $\alpha_{\dR}$ and
hence, by continuity, on $\omega \in\overline{\alpha}^p_{\dR}$.

Equation~\eqref{(2.2)} is a consequence of the minimality of $|\omega
|_p$ in $\overline{\alpha}^p_{\dR}$ and of the constancy of $|\omega
_x|$.  Note that we may assume that $\alpha\not= 0$ and hence $|\omega
_x|=r>0$ for almost all $x\in X$.  For every $\sigma\in\Omega^{k-1}X$
we have
$$
\begin{array}{rcl}
0=\left. \frac{d}{ds} \right|_{s=0} \left(\int_X |\omega +s\,
d\sigma|^p d \,\mbox{vol}\right) &=& p \int_X |\omega|^{p-2} \langle
\omega, d \sigma\rangle d \,\mbox{vol}\\ &=& p r^{p-2} \int_X \langle
\omega, d\sigma\rangle d \,\mbox{vol}.
\end{array}
$$
Since $r\not= 0$ this implies \eqref{(2.2)}.  Finally, we consider the
Euclidean bundle $\Lambda^* X = \bigoplus\limits^n_{k=0} \Lambda^k X$
and the selfadjoint, first order linear differential operator $d+d^*$
acting on sections of $\Lambda^* X$.  It is wellknown that the
operator~$d+d^*$ is elliptic.  Hence the standard regularity theory
for linear elliptic operators implies that $\omega$ is smooth and
$(d+d^*)(\omega)=0$ in the strong sense.  Since $(d+d^*)^2 = \Delta$,
we conclude that $\omega$ is harmonic.  From the fact that $\omega$ is
a harmonic form of constant norm, we can easily conclude that
$|\alpha|_q^*$ does not depend on $q$.  This is proved
in~\cite[Proposition~8.1]{BaKa03} for the case $k\in \{1, n-1 \}$, and
the proof carries over to the general case without any change.

Finally, we assume that $p\in [1,\infty)$, and that the form $\omega
\in \overline \alpha_{\dR}^p$ has constant norm almost everywhere and
satisfies $| \omega |_p = |\alpha |_p ^*$.  For every~$p' \in
(p,\infty]$, we have $\omega \in \overline \alpha_{\dR}^{p'}$ and
\[
| \alpha |_{p'}^* \leq | \omega | _{p'} = | \omega | _{p} = | \alpha
|_{p}^*  \leq | \alpha |_{p'}^* .
\]
Hence the first part of the proof shows that $| \alpha | _q^*$ does
not depend on $q$.  Since $ | \omega |_q$ is independent of $q$ as
well, we obtain
\[
| \alpha | _{p}^* = | \omega | _{p} = | \omega |_{q} = | \alpha
|_{q}^* 
\]
for all $q\in [1,\infty]$.
\end{proof}

\section{The BI construction adapted to conformal norms}
\label{bi}

We now continue with the notation of Section~\ref{three}.  Consider
the norm~$\|\;\|_{\conf/\varphi}$ in $V$.  We apply Lemma~\ref
{ellipsoid-decomposition} to the norm $\|\;\|_{\conf/\varphi}$.  This
yields a decomposition
\[
\|\;\|_E^2 = \sum_{i=1}^N \lambda_i L_i^2,
\]
where $\lambda_i>0$, $\sum \lambda_i=n$, $L_i\in V^*$ and
$\|L_i\|^*_{\conf/\varphi}=1$.  Then a linear map $L:V\to\R^N$ defined
by
\begin{equation}
\label{41b}
L(x) = (\lambda_1^{1/2} L_1(x), \lambda_2^{1/2} L_2(x), \dots,
\lambda_N^{1/2} L_N(x))
\end{equation}
is an isometry from $(V,\|\;\|_E)$ onto a subspace $L(V)$ of
$\R^N$, equipped with the restriction of the standard coordinate
metric of $\R^N$.

\begin{proposition}
\label{42b}
For every~$L_i$ there exists a weakly closed 1-form $\omega_i$ on $X$ from
the pullback cohomology class $\omega_i\in \varphi^*(L_i)$ and
\[
\|\omega_i\|_{n}=\|L_i\|^*_{\conf/\varphi}=1.
\]
\end{proposition}

\begin{proof}
The proof results from Proposition~\ref{Proposition1.3.}.  Briefly,
the issue is that the minimizing form $\omega_i$ will {\em a priori\/}
only be an $L^n$-form, \ie the $L^n$-limit of a sequence of smooth
forms in the pullback class $\varphi^*(L_i)$.  Then, it is a
non-trivial result from nonlinear PDE that $\omega_i$ is H\"{o}lder.
This can be found in \cite{To}, and in \cite{Ha92} for the case of
forms of arbitrary degree.  
\end{proof}

\begin{definition} 
\label{xx}
Denote by $\XX$ the covering space of $X$ defined by the subgroup
$ker(\varphi_*)\subset\pi_1(X)$.  We denote the action of $\Gamma =
\pi_1(X)/\ker(\varphi_*)$ on $\XX$ by $(v,x)\mapsto x + v$.
\end{definition}

\begin{corollary}
There exists a $C^1$ function
\begin{equation}
\label{21}
f_i:\XX\to\R \hbox{\ such\ that\ } df_i=\tilde\omega_i ,
\end{equation}
where $\tilde\omega_i$ is the lift of $\omega_i$.  Furthermore, $f_i$
satisfies the relation
\begin{equation}
\label{41}
f_i(x+v)=f_i(x)+L_i(v)
\end{equation}
for all $x\in\XX$ and $v\in\Gamma$. 
\end{corollary}

\begin{proof}
The proof of Proposition~\ref{42b} implies the existence of $f_i$ with
H\"{o}lder continuous first derivatives satisfying \eqref{21}.  Then
\eqref{41} follows from the fact that $\omega_i\in \varphi^*(L_i)$.
\end{proof}

Now we use the functions $f_i:\XX\to\R$ of formula \eqref{21} to
define a map~$F:\XX\to\R^N$ as follows:
\begin{equation}
\label{42}
F(x) = (\lambda_1^{1/2} f_1(x), \lambda_2^{1/2} f_2(x), \dots,
\lambda_N^{1/2} f_N(x)) .
\end{equation}
Observe that both the map $L$ of \eqref{41b} and $F$ of \eqref{42} are
$\Gamma$-equivariant with respect to the following action of $\Gamma$
on $\R^N$:
\[
\Gamma\times\R^N \to \R^N, \quad (v,x) \mapsto x+L(v).
\]
Now let $\Pr_{L(V)}:\R^N \to L(V)$ be the orthogonal projection to the
image~$L(V)$.  Then the composition 
\begin{equation}
L^{-1}\circ \Pr_{L(V)}\circ F : \XX \to V
\end{equation}
is a $\Gamma$-equivariant map covering a map $f: X \to V/\Gamma$.  We
will exploit it in Section~\ref{five}.

\section{Proof of Proposition~\ref{21a}}
\label{five}

The proof is similar to the one in \cite{IK}, with some crucial
differences.  Therefore we present the modified proof here.  The main
difference is that the functions~$f_i$ in formula \eqref{21} were
defined differently in \cite{IK}.  Consider the norm $\|\;\| _{\conf/
\varphi}$.  We apply the construction of Section~\ref{bi} to this
norm.  Since the projection is nonexpanding and the map $L$ is an
isometry, it suffices to prove the proposition for the map~$F$ of
formula~\eqref{42}.

Let $A=dF:T_x X \to \R^N$.  Then
$$
\trace(A^*A) = \sum \lambda_i |df_i|^2 .
$$
By the inequality of geometric and arithmetic means, we have
\[
\begin{aligned}
\jac(F) (x) &= \det(A^*A)^{1/2} \\ & \le \left(\tfrac 1n
\trace(A^*A)\right)^{n/2} .
\end{aligned}
\]

Thus, the map $F:\XX\to\R^N$ satisfies
\begin{equation}
\label{51b}
\jac(F) \le \left(\tfrac 1 n \sum \lambda_i |df_i|^2\right)^{n/2} .
\end{equation}
Therefore,
\[
\begin{aligned}
\jac(F) & \leq \left( \sum \tfrac{\lambda_i}{n}
|\tilde\omega_i|^2\right) ^{n/2} & \\ & \le \sum \tfrac{\lambda_i}{n}
|\tilde\omega_i|^n .
\end{aligned}
\]
The last inequality follows from Jensen's inequality applied to the
function $t\mapsto t^{n/2}$, when $n\ge 2$.  In the case $n=1$, there
is only one term in the ellipsoid decomposition, so Jensen is not
needed.

We now integrate the inequality over $X$.  Note that $\int_X
|\omega_i|^n=1$.  This yields the desired inequality for the area of
the corresponding map $X\to V/\Gamma$.

\section{Proof of Theorem \ref{13a} and  Theorem~\ref{35}}
\label{six}

\begin{proof}[Proof of Theorem~\ref{13a}]

Note that equality in \eqref{51b} is possible only if $F$ is
conformal.  The discussion of equality in \eqref{51b} leads to the
fact that equality in \eqref{10c} implies that the map $f$ is
conformal at all points where the differential does not vanish.  If
$n\geq 3$, then by \cite{Fer}, the map $f$ is indeed a smooth local
(conformal) diffeomorphism, hence a covering.  Here we are assuming
that the metric on $X$ is smooth.  Hence~$X$ is indeed a torus and $f$
is a conformal diffeomorphism.

If~$n=2$, then we can lift $f$ to a holomorphic map between the
universal covers.  Since $f$ is not constant, we see that both have to
be the complex plane.  Since the Klein bottle does not qualify as $X$,
we see that $X$ is a 2-torus, while $f$ is a nonconstant holomorphic
map between 2-tori.  So its lift to the universal covers is affine,
and hence $f$ is a covering in this case, as well.
\end{proof}

\begin{proof}[Proof of Theorem~\ref{35}]
We apply Proposition~\ref{125} together with Theorem~\ref{11c} (which
is immediate from Proposition~\ref{21a}) and an argument based on the
perfection of a critical lattice similar to the proof of~\cite
[Proposition~10.5] {BK}.  In more detail, let us normalize the metric
to unit volume to fix ideas.  The equality of the stable norm and the
conformally invariant norm for each of the integral elements
realizing~$\lambda_1 \left( H_1(X,\Z)_\R, \|\;\| \right)$ proves that
the harmonic 1-form in each of those integral classes has pointwise
constant norm.  Next, perfection implies that in fact {\em every\/}
harmonic 1-form on $X$ has pointwise constant norm.  It follows that
the standard Abel-Jacobi map defined by those forms, is a Riemannian
submersion onto a flat torus, without any ramification.  Since the
induced homomorphism of fundamental groups is surjective, the degree
of the submersion is 1, proving the theorem.
\end{proof}

\section{Case $n \geq b$ and $L^p$ norms in homology}

In this section, we generalize the techniques of Section~\ref{bi}.

\begin{definition}
Let $X^n$ and $Y^\dimtarget$ be Riemannian manifolds, $n\ge
\dimtarget$, and $f\colon X\to Y$ a smooth map. For every $x\in X$, we
consider the derivative~$d_x f$, and define
$$
\jac^\perp f(x) = \sup_\sigma \frac {\vol_\dimtarget (d_xf(\sigma))}
{\vol_\dimtarget (\sigma)}
$$
where $\sigma$ ranges over all $\dimtarget$-dimensional parallelotopes
in $T_x X$.
\end{definition}

Obviously, if $x$ is a regular point of $f$, then the above supremum
is attained at $\sigma$ lying in the orthogonal complement of the
fiber $f^{-1}(f(x))$, and $\jac^\perp f(x)=0$ if $x$ is not a regular
point.  We will use the following coarea formula.  Every smooth map
$f\colon X^n\to Y^ \dimtarget$ satisfies the identity
\begin{equation}
\label{51}
\int_X \jac^\perp f(x) \; d \vol_n(x) = \int_Y \vol_{n-\dimtarget}
(f^{-1}(y)) \ d\vol_\dimtarget (y) .
\end{equation}

We present a version of Proposition \ref{area-nonexp} and
Proposition~\ref{21a}, which allows us to prove inequality~\eqref{23c}
for the conformal systole, as well as to handle the case of equality
when $n \geq b$.

\begin{proposition} 
\label{area-nonexp-int}
Let $X^n$ be a closed Riemannian manifold of unit volume, and $V ^
\dimtarget$ a vector space, $n\ge \dimtarget$.  Let $p \geq
\max\{\dimtarget ,2 \}$.  Let $\Gamma$ be lattice in $V$, and let
$\varphi:X\to V/\Gamma$ be a continuous map inducing an epimorphism of
the fundamental groups.  Consider the associated~$L^p$-norm $\|\;\|_p$
on $V$, \cf Definition~\ref{23z}.  Denote by~$\| \; \| _E$ the
Euclidean norm on $V$, defined by the John ellipsoid
of~$\|\;\|_{\stp}$.  Then the homotopy class of $\varphi$ contains a
smooth map $f:X\to(V/\Gamma,\|\;\|_E)$ such that
\begin{equation}
\label{f1}
\int_X \jac^\perp f \le 1 .
\end{equation}
\end{proposition}

\begin{proposition}
\label{103}
In the hypotheses of Proposition~\ref{area-nonexp-int}, assume a
strict inequality $p>\max\{\dimtarget ,2 \}$.  If the equality is
attained in \eqref{f1}, then the map $f$ is a harmonic Riemannian
submersion.
\end{proposition}
In section~\ref{eight}, we will generalize the BI construction of
Section~\ref{bi}, and use it in Section~\ref{nine} to prove the
Propositions~\ref{area-nonexp-int} and \ref{103}.

\section{The BI construction in the case $n\geq b$} 
\label{eight}

Applying Lemma \ref{ellipsoid-decomposition} to the norm $\| \; \|
_{\stp}$ yields a decomposition
\[
\|\;\|_E^2 = \sum_{i=1}^N \lambda_i L_i^2,
\]
where $\lambda_i>0$, $\sum \lambda_i=\dimtarget$, $L_i\in V^*$ and
$\|L_i\|^*_{\stp}=1$.  Then a linear map~$L:V\to\R^N$ defined by
$$
L(x) = (\lambda_1^{1/2} L_1(x), \lambda_2^{1/2} L_2(x), \dots,
\lambda_N^{1/2} L_N(x))
$$
is an isometry from $(V,\|\;\|_E)$ onto a subspace $L(V)$ of
$\R^N$, equipped with the restriction of the standard coordinate
metric of $\R^N$.

We identify the dual space $V^*$ of $V$ with $H^1(V/\Gamma;\R)$.
Then, we have a pull-back map $\varphi^*:V^*\to H^1(X;\R)$.  For every
$\alpha\in V^*$, we have
$$
\|\alpha\|^*_{\stp} = \|\varphi^*(\alpha)\|^*_{\stp} = \inf \left\{
|\omega|_p \left|\; \omega\in\varphi^*(\alpha)_{\dR}^{\phantom{I}}
\right. \right\} .
$$
\begin{definition} 
\label{53}
Fix a real $p> \max\{\dimtarget,2\}$ and let $\omega_i$, $i= 1,
\ldots, N$ denote the $L^p$-minimizer in the cohomology class
$\varphi^* (L_i) \in H^1(X;\R)$.
\end{definition}

Since $\vol(X)=1$, we have $|\omega_i|_p=\|L_i\|^*_{\stp}=1$.  The rest of
the construction is the same as in Section~\ref{bi}.  In particular,
we define a map
\begin{equation}
\label{81}
F(x) = (\lambda_1^{1/2} f_1(x), \lambda_2^{1/2} f_2(x), \dots,
\lambda_N^{1/2} f_N(x)) .
\end{equation}
In conclusion,
we obtain a map $\tilde f\colon\XX\to V$ defined by
\begin{equation}
\label{f5}
\tilde f=L^{-1}\circ \Pr_{L(V)}\circ F
\end{equation}
where $\Pr_{L(V)}:\R^N \to L(V)$ is the orthogonal projection
to~$L(V)$.  Clearly $\tilde f$ is $\Gamma$-equivariant, hence it is a
lift of a map $f\colon X\to V/\Gamma$.  We show in Section~\ref{nine}
that $f$ satisfies inequality \eqref{f1}.

\section{Proof of Propositions~\ref{area-nonexp-int} and \ref{103}}
\label{nine}
Consider the map $f$ constructed in Section~\ref{eight}.  Let $x\in
X$.  If $x$ is a regular point of $f$, we let $\Sigma_x=(\ker
d_xf)^\perp$.  Otherwise let $\Sigma_x$ be an arbitrary
$d$-dimensional linear subspace of $T_x X$.  Then
\begin{equation}
\label{f6}
\jac^\perp f(x) = \jac ((d_x f)|_{\Sigma_x}) \le \jac ((d_x
F)|_{\Sigma_x})
\end{equation}
since the projection in \eqref{f5} is nonexpanding and $L$ is an
isometry.  Let~$A$ be the linear map
\begin{equation}
\label{47}
A=(d_x F)|_{\Sigma_x}:\Sigma_x\to\R^N .
\end{equation}
Then, by \eqref{81},
$$
\trace(A^*A) = \sum \lambda_i \left| (d_xf_i)|_{\Sigma_x}
^{\phantom{I}} \right|^2 \le \sum \lambda_i |d_xf_i|^2 = \sum
\lambda_i |\omega_i(x)|^2
$$
By the inequality of geometric and arithmetic means, we have
\begin{equation}
\label{f7}
\jac ((d_x F)|_{\Sigma_x}) = \det(A^*A)^{1/2} \le \left(\tfrac 1
\dimtarget \trace(A^*A)\right)^{\dimtarget /2} .
\end{equation}
Thus, the map $f$ satisfies
\begin{equation}
\label{f8}
\jac^\perp f \le \left(\sum\tfrac{\lambda_i}{\dimtarget} | \omega_i |
^2 \right)^{\dimtarget /2} \le \left(\sum \tfrac {\lambda_i}
{\dimtarget} |\omega_i|^p\right)^{\dimtarget/p}.
\end{equation}
The last inequality follows from Jensen's inequality applied to the
function $t\mapsto t^{p/2}$; recall that $p>2$ and
$\sum\lambda_i=\dimtarget$.  Integrating over $X$ yields
\begin{equation}
\label{510}
\int_X \jac^\perp f \le \int_X \left(\sum \tfrac{\lambda_i}
{\dimtarget} |\omega_i|^p\right)^{\dimtarget/p} \le \left(\int_X \sum
\tfrac{\lambda_i}{\dimtarget} |\omega_i|^p\right)^{\dimtarget/p}
\end{equation}
since $\dimtarget/p<1$ and $\vol(X)=1$.  Since $|\omega_i|_p\le 1$, we have
\begin{equation}
\label{f9}
 \int_X \sum \tfrac{\lambda_i}{\dimtarget} |\omega_i|^p
 = \sum \tfrac{\lambda_i}{\dimtarget} \int_X |\omega_i|^p
 \le \sum \tfrac{\lambda_i}{\dimtarget} = 1,
\end{equation}
and the desired inequality \eqref{f1} follows.


\begin{remark}
\label{121}
Let $f: X \to J_1$ be Lichnerowicz's harmonic map.  Let~$p=b=2$.
Let~$\omega_1, \omega_2$ be an orthonormal basis for the harmonic
1-forms with respect to the $L^2$ inner product.  
Then one has $\lambda_1=\lambda_2=1$, and the computations above
can be simplified  as follows:
\[
\int_X \jac^\perp f = \int_X |\omega_1\wedge\omega_2| \le \tfrac12
\int_X |\omega_1|^2+|\omega_2|^2=1,
\]
and the desired inequality \eqref{f1} follows.
\end{remark}

In the case of equality in \eqref{f1}, observe that all functions
under integrals are continuous, hence all the inequalities throughout
the argument turn to equalities.

\begin{lemma}
\label{91}
Equality in \eqref{f1} implies that $|\omega_i|=1$ everywhere.
\end{lemma}

\begin{proof}
Equality in \eqref{f8} implies that $|\omega_i|=|\omega_j|$ for
all $i,j$. The equality in \eqref{510} implies that the function
\[
\sum_i \frac{\lambda_i}{\dimtarget} |\omega_i|^p
\]
is constant. These together imply that $|\omega_i|$ is constant.
Then equality in \eqref{f9} implies that $|\omega_i|=1$.
\end{proof}

\begin{proof}[Proof of Proposition~\ref{103}]
Since, by Lemma~\ref{91}, we have $|\omega_i|=1$ everywhere and
$\omega_i$ is an $L^p$-minimizer, Proposition~\ref{125} implies
that~$\omega_i$ is harmonic for every $i$.  Hence $f$ is harmonic
(and, in particular, smooth).

Equality in \eqref{f7} is attained precisely when the map $A$ of
formula~\eqref{47} is a conformal linear map.  Then equality in
\eqref{f6} implies that $(df)|_{\Sigma_x}$ is conformal, too.
Finally, equality in \eqref{f8} implies that~$\jac^\perp f=1$, thus
$f$ is a Riemannian submersion.
\end{proof}

\section{Proof of Proposition~\ref{23b} and Theorem \ref{17}}
\label{thirteen}

\begin{proof}[Proof of Proposition \ref{23b}]
We apply Proposition \ref{area-nonexp-int} to the Abel-Jacobi
map~$\AJ:X\to \T^b=H_1(X;\R)/H_1(X;\Z)_\R$.  We obtain a map $f$ which
is homotopic to $\AJ$ and satisfies inequality \eqref{f1}, where
$\T^b$ is equipped with the flat Euclidean metric defined by the John
ellipsoid of the unit ball of the the norm $\|\;\|_{p}$, where the
case $p=\infty$ corresponds to the stable norm.  Then, by the coarea
formula \eqref{51},
$$
\int_X \jac^\perp f(x) \; d \vol_n(x) = \int_{\T^b}
\vol_{n-b}(f^{-1}(y)) \ d\vol_b(y) .
$$
Observe that if $y\in \T^b$ is a regular value for $f$, then
\begin{equation}
\label{fiber-degree}
\vol_{n-b}(f^{-1}(y))\ge\deg(\AJ),
\end{equation}
as in Definition \ref{deg}.  Hence we have
$$
\deg(\AJ) \vol_b \left( \T^b, \|\;\|_E \right) \le \vol_n(X).
$$
by inequality \eqref{f1}.  The proposition now results from the
definition of the Hermite constant in formula \eqref{34}.
\end{proof}

\begin{proof}[Proof of Theorem \ref{17}]
The theorem is immediate from Proposition~\ref{23b} combined with
Proposition~\ref{125}, by the argument of the proof of
Theorem~\ref{35} in Section \ref{six}.  Alternatively, one could argue
as follows.  In the case of equality, the equalities in \eqref{f1} and
\eqref{fiber-degree} are attained.  This implies the inequality
\eqref{eq12}, and the theorem follows from Lemma~\ref{91}
characterizing the equality case in Proposition~\ref{area-nonexp-int}.
Note that equality in \eqref{fiber-degree} implies that $f^{-1}(y)$ is
minimal.
\end{proof}

\section{Acknowledgments} 
We are grateful to B. White for a helpful discussion of the material
of Section~\ref{smallcycles}.

\vfill\eject

\end{document}